# Local piecewise polynomial curve of variable order and its appliance for calculating of nodal derivatives

Oleg Stelia[a,*], Leonid Potapenko[a], Igor Stelia[b]

[a] Faculty of Computer Science and Cybernetics, Taras Shevchenko National University of Kyiv, 03680, Ukraine
[b] EPAM Systems, Kyiv, 03150, Ukraine

**Abstract**

In this paper, we present a new Hermite type curve piecewise polynomial of $C^1$ continuity on nonuniform grids. Depending on the location of the knots of the grid (uniform and nonuniform), the curve is quadratic or cubic. The approximation error is investigated. For nonuniform grids a method for calculating the first and second derivatives in knots is proposed. The order of approximation of the first derivative is $O(h^2)$ and of the second derivative is $O(h)$. We also perform several numerical experiments which exemplify the properties of the proposed curve. Graphic examples are given.

## 1. Introduction

Let $\Delta_\tau : a = \tau_1 < \tau_2 < ... < \tau_N = b$ be a partition of an interval $[a,b]$. The problem of constructing the cubic Hermite interpolant [1] is as follows: it is necessary to construct polynomials $s_i(x)$, $x \in [\tau_i, \tau_{i+1}]$ at intervals $[\tau_i, \tau_{i+1}]$, $i = \overline{1, N-1}$, for which the following conditions are satisfied at the knots of the grid $\Delta_\tau$:

$$s_i(\tau_i) = F_i, \; s_i(\tau_{i+1}) = F_{i+1}, \; s_i'(\tau_i) = F_i', \; s_i'(\tau_{i+1}) = F_{i+1}', \; i = \overline{1, N-1}, \quad (1)$$

where the values $F_i$ are given data, the derivatives $F_i'$ are free parameters.

The variety of methods proposed in the literature for constructing such polynomials is caused by the various requirements that apply to them. Such requirements include: the degree of polynomials, a piecewise polynomial curve is local or global, the requirement for a smoothness of a curve, monotonicity and / or convexity, and others [2].

Algorithms for preserving the interpolation shape using quadratic splines were proposed in [3], [4]. In [3], a quadratic piecewise polynomial is constructed using Bernstein polynomials. A similar algorithm with quadratic piecewise polynomials is presented in [4]. This algorithm preserves monotonicity and / or convexity by adding one node, if necessary, in each data subinterval. In [5], an algorithm is described in which the initial piecewise cubic interpolant is modified by changing the derivative values of the interpolant (if necessary) to produce a monotonic piecewise cubic interpolant. In [6] and [7], the derivative values are not modified. To ensure the monotonicity, additional knots are inserted into the subinterval. The Fritsch and Carlson method (1980) [5] was developed in [8], where a new technique for slope estimation was proposed. In [9], a method is described for producing monotone piecewise cubic interpolants for monotone data that are completely local. An algorithm that, if necessary, adds one or two

*Corresponding author.
E-mail addresses: oleg.stelya@gmail.com (O. Stelia), lpotapenko@ukr.net (L. Potapenko), igor.stelia@gmail.com (I. Stelia)
This research did not receive any specific grant from funding agencies in the public, commercial, or not-for-profit sector additional knots between existing knots in order to preserve the monotonicity is presented in [10]. A new optimal property, namely, the minimum quadratic oscillation on average was obtained in [11] for the Hermitian type cubic splines. The method requires solving a system with a dominant diagonal matrix. To preserve the monotonicity an interpolant was developed in [12] on the rational Bezier cubic basis in which sufficient conditions are expressed in terms of weights. New methods for assigning derivatives at discontinuity points are presented and analyzed in [13] and [14], which lead to third-order approximations even in the case of nonuniform grids. Reviews of various methods and calculation examples are given in [15].

In this paper, we propose a new $C^1$ continuity local piecewise polynomial curve of variable degree. In the general case, the curve does not pass through the given points. A method for calculating the first and second derivatives at the grid knots is proposed. We describe the problem in Section 2, discuss methods in Section 3, and present test results in Section 4.

## 2. Construction of the proposed curve

In this section, we present an algorithm for the construction of curves approximating a set of data points.

Let us define the grid $\Delta_x : \tau_1 \leq x_2 < ... < x_N \leq \tau_N$. We will construct polynomials $s_i(x)$, $i = \overline{2, N}$, $x \in [x_i, x_{i+1}]$, $i = \overline{2, N-1}$, satisfying the conditions:

$$s_i(x_i) = f_i, \ s_i(x_{i+1}) = f_{i+1}, \ s_i'(x_i) = f_i', \ s_i'(x_{i+1}) = f_{i+1}', \tag{2}$$

where

$$x_i = \alpha_i \tau_{i-1} + (1 - \alpha_i) \tau_i, \ 0 < \alpha_i < 1, \ i = \overline{3, N-1}, \ 0 < \alpha_2 \leq 1,$$
$$x_{i+1} = \beta_i \tau_{i+1} + (1 - \beta_i) \tau_i, \ 0 < \beta_i < 1, \ \alpha_{i+1} = 1 - \beta_i, \ i = \overline{2, N-2}, \ 0 < \beta_{N-1} \leq 1,$$
$$f_i = \alpha_i F_{i-1} + (1 - \alpha_i) F_i,$$
$$f_{i+1} = \beta_i F_{i+1} + (1 - \beta_i) F_i, \tag{3}$$
$$f_i' = (F_i - F_{i-1}) / H_i, \text{ where } H_i = \tau_i - \tau_{i-1},$$
$$f_{i+1}' = (F_{i+1} - F_i) / H_{i+1}, \text{ where } H_{i+1} = \tau_{i+1} - \tau_i. \tag{4}$$

Then the polynomial $s_i(x)$ on the interval $[x_i, x_{i+1}]$ is written as follows:

$$s_i(x) = f_{i+1} \frac{x - x_i}{x_{i+1} - x_i} + f_i \frac{x - x_{i+1}}{x_i - x_{i+1}} + (x - x_i)(x - x_{i+1})(Ax + B), \tag{5}$$

where

$$A = \frac{f_i' + f_{i+1}'}{(h_i + h_{i+1})^2} - \frac{2(f_{i+1} - f_i)}{(h_i + h_{i+1})^3},$$

$$B = \frac{f_{i+1}' - f_i'}{2(h_i + h_{i+1})} - \frac{(f_{i+1}' + f_i')(x_i + x_{i+1})}{2(h_i + h_{i+1})^2} + \frac{(f_{i+1} - f_i)(x_i + x_{i+1})}{(h_i + h_{i+1})^3},$$

$$h_i = \alpha_i H_i, \ h_{i+1} = \beta_i H_{i+1}.$$

**Remark**. For $h_i = h_{i+1}$ and $H_i = H_{i+1}$ we have $A = 0$, then therefore the polynomials $s_i(x)$ will be quadratic.

A piecewise cubic polynomial curve $S(x) \in C^1[a, b]$ formed by polynomials $s_i(x)$, $i = \overline{2, N}$, in the general case, does not pass through the points $F_i$. If you put $x_2 = \tau_1$ and $x_N = \tau_N$ the curve will start at a point $F_1$ and end at a point $F_N$.

In contrast to the global curve presented in [16], for constructing a local curve there is no need to solve a system of equations and there are no restrictions on the $x_i$ knots location.

## 3. Approximation error

In this section, we consider the approximation error. Let
$$\delta_1 = S(\tau_i) - F(\tau_i), \quad \delta_2 = S'(\tau_i) - F'(\tau_i).$$
We assume that $F(x)$ is sufficiently smooth on a interval $[a,b]$. Using the Taylor formula, taking into account expressions (3) - (5), we obtain for uniform grid ($h_i = h_{i+1} = h$, $H_i = H_{i+1}$):

$$\delta_1 = \frac{1}{2}h^2 F'(\tau_i) + O(h^4), \tag{6}$$

$$\delta_2 = \frac{2}{3}h^2 F'''(\tau_i) + O(h^3), \tag{7}$$

$$S''(\tau_i) - F''(\tau_i) = \frac{2}{3}h^2 F^{(4)}(\tau_i) + O(h^4).$$

For a uniform grid, $S'(\tau_i)$ and $S''(\tau_i)$ coincide with the first and second derivatives of the quadratic function passing through the points $(\tau_{i-1}, F_{i-1}), (\tau_i, F_i), (\tau_{i+1}, F_{i+1})$.

For nonuniform grid ($h_i \neq h_{i+1}$, $H_i \neq H_{i+1}$, $\bar{h} = \max(h_i, h_{i+1})$).

$$\delta_1 = C_1 F''(\tau_i) + O(\bar{h}^3), \tag{8}$$

where

$$C_1 = h_i h_{i+1}[(\frac{h_{i+1}}{\beta} + \frac{h_i}{\alpha})/4(h_i + h_{i+1}) + (\frac{h_{i+1}}{\beta} - \frac{h_i}{\alpha})(h_{i+1} - h_i)/4(h_i + h_{i+1})^2$$

$$-(\frac{h_{i+1}^2}{\beta} - \frac{h_i^2}{\alpha})(h_{i+1} - h_i)/2(h_i + h_{i+1})^3].$$

$$\delta_2 = C_2 F''(\tau_i) + O(\bar{h}^2), \tag{9}$$

$$C_2 = (\frac{h_{i+1}^2}{\beta} - \frac{h_i^2}{\alpha})/2(h_i + h_{i+1}) + (h_{i+1}^2 + h_i^2 - 4h_i h_{i+1})(\frac{h_{i+1}}{\beta} - \frac{h_i}{\alpha})/4(h_i + h_{i+1})^2 -$$

$$(h_{i+1}^2 + h_i^2 - 4h_i h_{i+1})(\frac{h_{i+1}^2}{\beta} - \frac{h_i^2}{\alpha})/2(h_i + h_{i+1})^3 - (h_{i+1} - h_i)(\frac{h_{i+1}}{\beta} + \frac{h_i}{\alpha})/4(h_i + h_{i+1}).$$

To calculate the derivatives at the $\Delta_\tau$ grid knots using formulas (8) and (9) we obtain:

$$F''(\tau_i) = (S(\tau_i) - F(\tau_i))/C_1 + O(\bar{h}), \tag{10}$$

$$F'(\tau_i) = S'(\tau_i) - C_2 F''(\tau_i) + O(\bar{h}^2). \tag{11}$$

We formulate the results in this section as a lemma.

**Lemma**. Let $F(x)$ – the function, sufficiently smooth along with its derivatives on the interval $[a,b]$, $S(x)$ - be a constructed piecewise polynomial curve. Then the following inequalities are performed:

$$|S(\tau_i) - F(\tau_i)| \leq \begin{cases} O(h^2), \text{ for the uniform grid}, \\ O(\bar{h}^2), \text{ for the nonuniform grid}. \end{cases}$$

$$|S'(\tau_i) - F'(\tau_i)| \leq \begin{cases} O(h^2), \text{ for the uniform grid}, \\ O(\bar{h}^2), \text{ for the nonuniform grid}. \end{cases}$$

$$|S''(\tau_i) - F''(\tau_i)| \leq \begin{cases} O(h^2), \text{ for the uniform grid}, \\ O(\bar{h}), \text{ for the nonuniform grid}. \end{cases}$$

## 4. Numerical experiments

In this section we will perform some examples using the formulas obtained in the previous section.

**Example 1.** We use an example of a smooth function $F(x) = x^4 + \sin(x)$, let $\tau_{i-1} = 0.5 - H_i$, $\tau_i = 0.5$, $\tau_{i+1} = 0.5 + H_{i+1}$, where $H_i = 2^{-j}$, $5 \leq j \leq 9$ [13].
The errors are calculated as

$$err_1 = |S(\tau_i) - F(\tau_i)|,$$
$$err_2 = |S'(\tau_i) - F'(\tau_i)|,$$
$$err_3 = |S''(\tau_i) - F''(\tau_i)|.$$

The estimated order of convergence (EOC) is calculated via $\log_2(err(2^{-(j-1)})/err(2^{-j}))$.
Results for uniform and nonuniform grids are presented in Table 1 and Table 2 respectively.

Table 1

Errors and estimated orders of convergence. Uniform grid ($H_{i+1} = H_i$, $\alpha_i = \beta_i = 1/2$).

| $H_i$ | Results for $F(\tau_i)$ | | Results for $F'(\tau_i)$ | | Results for $F''(\tau_i)$ | |
|---|---|---|---|---|---|---|
| | $err_1$ | EOC | $err_2$ | EOC | $err_3$ | EOC |
| 3.125e-2 | 3.0793e-4 | - | 1.8102e-3 | - | 1.9921e-3 | - |
| 1.562e-2 | 7.6937e-5 | 2.0009 | 4.5257e-4 | 1.9999 | 4.9803e-4 | 2.0000 |
| 7.812e-3 | 1.9231e-5 | 2.0002 | 1.1314e-4 | 2.0000 | 1.2450e-4 | 2.0001 |
| 3.906e-3 | 4.8076e-6 | 2.0000 | 2.8285e-5 | 2.0000 | 3.1127e-5 | 1.9999 |
| 1.953e-3 | 1.2019e-6 | 2.0000 | 7.0714e-6 | 2.0000 | 7.7819e-6 | 2.0000 |

Table 2

Errors and estimated orders of convergence. Nonuniform grid ($H_{i+1} = 3H_i$, $\alpha_i = \beta_i = 1/2$).

| $H_i$ | Results for $F(\tau_i)$ | | Results for $F'(\tau_i)$ | | Results for $F''(\tau_i)$ | |
|---|---|---|---|---|---|---|
| | $err_1$ | EOC | $err_2$ | EOC | $err_3$ | EOC |
| 3.125e-2 | 7.5977e-4 | - | 5.6177e-3 | - | 2.4567e-1 | - |
| 1.562e-2 | 1.8126e-4 | 2.0675 | 1.3810e-3 | 2.0243 | 1.1934e-1 | 1.0416 |
| 7.812e-3 | 4.4277e-5 | 2.0334 | 3.4234e-4 | 2.0122 | 5.8800e-2 | 1.0212 |
| 3.906e-3 | 1.0942e-5 | 2.0167 | 8.5222e-5 | 2.0061 | 2.9182e-2 | 1.0107 |
| 1.953e-3 | 2.7198e-6 | 2.0083 | 2.1259e-5 | 2.0031 | 1.4536e-2 | 1.0054 |

**Example 2.** Now we give an example of graphs produced by the Algorithms of Section 2 and Section 3. We use data set considered by [5].

Table 3
Data for Example 2.

| $i$ | 1 | 2 | 3 | 4 | 5 | 6 | 7 | 8 | 9 |
|---|---|---|---|---|---|---|---|---|---|
| $\tau_i$ | 7.99 | 8.09 | 8.19 | 8.7 | 9.2 | 10.0 | 12.0 | 15.0 | 20.0 |
| $F_i$ | 0.0 | 0.0000276429 | 0.0437498 | 0.169183 | 0.469428 | 0.94374 | 0.998636 | 0.999919 | 0.999994 |

Depending on the $x_i$ point location, we obtain various curves approximating the initial data (Fig. 1.).

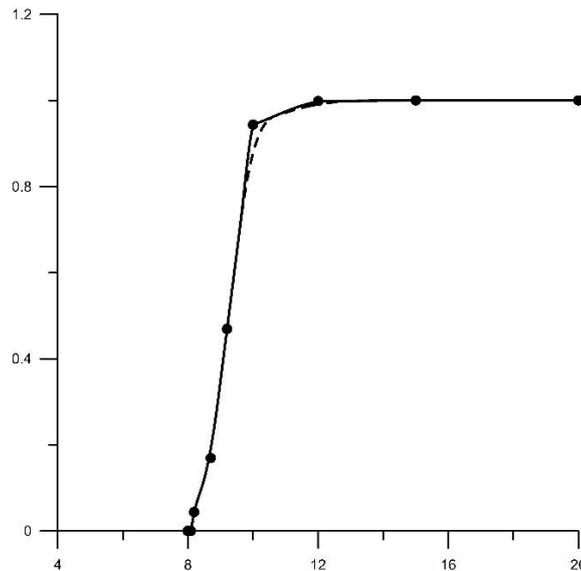

Fig. 1. Reconstructions from data of Table 3.
Experiment 1 (dotted line), Experiment 2 (solid line).

Table 4
$x_i$ knots location.

| | $i$ | 2 | 3 | 4 | 5 | 6 | 7 | 8 | 9 |
|---|---|---|---|---|---|---|---|---|---|
| Exp. 1 | $x_i$ | 7.99 | 8.14 | 8.445 | 8.95 | 9.6 | 11.0 | 13.5 | 20.0 |
| Exp. 2 | $x_i$ | 7.99 | 8.14 | 8.445 | 8.95 | 9.6 | 10.1 | 12.1 | 20.0 |

## 5. Conclusion

Using the polynomials forming the proposed curve, we can calculate the first and second derivatives of the grid function. In practice, derivatives are usually most interesting at the same points where the function itself is specified. Control of the knots $x_i$ provides additional freedom in constructing a curve. As a result, various curves are obtained, among which the most appropriate curve is selected.

## References

[1] C. de Boor, A Practical Guide to Splines, Springer-Verlag, Berlin, 2001.